\documentclass[11pt,openbib]{article}

\usepackage[T1]{fontenc}
\usepackage{amsmath}
\usepackage{amsfonts}
\usepackage{amssymb}
\usepackage{graphicx}

\newcommand{\R}{\mathbb{R}} 
\newcommand{\C}{\mathbb{C}}

\newcommand{\Z}{\mathbb{Z}}

\newcommand{\spp}{\mathop{\mathrm{sp}}\nolimits}

\newcommand{\dee}{\mathop{\! \,\mathrm{d} \!}\nolimits}
\newcommand{\comp}{\raisebox{0pt}{$\scriptstyle\circ \, $}}
\newcommand{\setrule}{\, \rule[-4pt]{.5pt}{13pt}\, }

\newcommand{\smallspace}{\smallskip\par\noindent}

\newcommand{\onehalf}{\mbox{$\frac{\scriptstyle 1}{\scriptstyle 2}\,$}} 
 
\newcommand{\lefthook}{\mbox{$\, \rule{8pt}{.5pt}\rule{.5pt}{6pt}\, \, $}}
\newcommand{\ttfrac}[2]{\mbox{$\frac{{\scriptstyle #1}}{{\scriptstyle #2}}$}}

\newcommand{\vvee}{\mbox{\tiny $\vee $}}

\allowdisplaybreaks

\begin{document}
\begin{center}
{\Large \bf Geometric scattering monodromy} \\
\vspace{.05in} 
Richard Cushman\footnotemark
\end{center}
\footnotetext{Department of Mathematics and Statistics, University of Calgary, Calgary, Alberta T2N 1N4, Canada. email: r.h.cushman@gmail.com; 
}
\addtocounter{footnote}{1}
\footnotetext{\today} \bigskip 

\noindent \textbf{Keywords} scattering monodromy, complex Morse lemma \medskip 

\noindent MSC (2020) 70H06 \medskip 

\begin{abstract}
In this paper we give geometric conditions so that the integral mapping of a Liouville integrable 
Hamiltonian system with a focus-focus equilibrium point has scattering monodromy. Using a 
complex version of the Morse lemma, we show that scattering monodromy is the same as 
the scattering monodromy of the standard focus-focus system.
\end{abstract} \bigskip 

\noindent {\Large \bf Introduction} 
\bigskip 

In \cite{bates-cushman} the hyperbolic oscillator integrable Hamiltonian system 
$(u,v, {\R }^4, \Omega $ $= \dee {\xi}_1 \wedge \dee {\xi}_2 - \dee {\eta }_1 \wedge \dee {\eta}_2)$, 
where 
\begin{align*}
&u: {\R }^4 \rightarrow \R : (\xi , \eta ) \mapsto {\xi }_1{\eta }_1+{\xi }_2{\eta }_2  = h \\
& \hspace{-1.1in} \mathrm{and} \\
& v: {\R }^4 \rightarrow \R : (\xi , \eta ) \mapsto \onehalf ( {\xi }^2_1 +{\xi }^2_2 - {\eta }^2_1 - {\eta }^2_2) = 
\ell 
\end{align*}
was shown to have scattering monodromy. Geometrically this means that a motion in ${\R }^4$ of the 
hyperbolic oscillator of energy $h$ and angular momentum $\ell $ projects onto a branch of a hyperbola 
in the $({\xi }_1, {\xi }_2)$ plane, whose outgoing asymptote forms an angle ${\tan }^{-1}\frac{h}{\ell }$ 
with its incoming asymptote. This angle is called the scattering angle of the hyperbolic motion. 
As $(h, \ell )$ traverses a circle in the energy-momentum plane centered at the origin, the scattering 
angle of a hyperbolic motion, which starts at a point in the image of a section of the bundle formed 
by the integral map, increases by $2\pi $. This is the scattering monodromy of 
the hyperbolic oscillator system. \medskip 

In \cite{dullin-waalkens} it was shown that the quantum Kepler problem has scattering 
monodromy. In \cite{efstathiou-giacobbe-mardesic-sugny} the relation of scattering monodromy 
to the geometric monodromy of a toral fibration was treated using rotation forms. However, 
a geometric scattering monodromy theorem was not formulated. This paper remedies this omission. \medskip 

Our formultation of the geometric scattering monodromy theorem follows that of the geometric 
(toral) monodromy theorem given in \cite{martynchuk-broer-efstathiou}. Our proof follows the line 
of argument for the proof of the toral geometric monodromy theorem given in \cite{cushman-bates} with 
all reasoning invoving compactness being avoided. We use a complex version of the Morse 
lemma, inspired by \cite{vungoc-wacheux}, to reduce the proof of the geometric scattering 
monodromy theorem to the computation of the scattering monodromy of the complexified 
standard focus-focus system $(q_1, q_2, {\R }^4, \omega = 
\dee x \wedge \dee p_x + \dee y \wedge \dee p_y)$, where 
\begin{align*}
&q_1: {\R }^4 \rightarrow \R : (x,y,p_x,p_y) \mapsto xp_x +yp_y \\
&\hspace{-1.1in}\mathrm{and} \\
& q_2: {\R }^4 \rightarrow \R :(x,y,p_x,p_y) \mapsto xp_y - yp_x. 
\end{align*}

\section{The geometric scattering monodromy theorem}

In this section we state the geometric scattering monodromy theorem. \medskip 

The origin $0$ of ${\R}^4$ is a \emph{focus-focus} equilibrium point of the Liouville 
integrable system $( h_1, h_2, {\R}^4, \omega = \dee x \wedge \dee p_x + \dee y \wedge \dee p_y)$ if and only if \medskip 

\indent 1. \parbox[t]{4.25in}{The complete vector fields $X_{h_1}$ and $X_{h_2}$ vanish at $0$, that is, 
$0$ is an equilibrium point of $X_{h_1}$ and $X_{h_2}$.} 
\smallspace
\indent 2. \parbox[t]{4.25in}{The space spanned by the linearized Hamiltonian vector fields 
$DX_{h_1}(0)$ and $DX_{h_2}(0)$ is conjugate by a real linear symplectic mapping of 
$({\R}^4, \omega )$ into itself to the Cartan subalgebra of $\spp (4, \R)$ spanned by $X_{q_1}$ and 
$X_{q_2}$, where $q_1 = xp_x + yp_y$ and $q_2 = xp_y -yp_x$.} \medskip 

\noindent From point 2 we may assume that $h_i = q_i +r_i$ for $i = 1,2$, where 
$r_i$ is a smooth function on ${\R}^4$, which is flat to $2^{\mathrm{nd}}$ order at $0$, that is, 
$r_i \in \mathcal{O}(2)$. \medskip 

\noindent The remainder of this paper is devoted to proving \medskip 

\noindent \textbf{Theorem (Geometric scattering monodromy)}. Let 
$( h_1, h_2, {\R}^4, \omega )$ be a Liouville integrable system with a focus-focus 
equilibrium point at $0 \in {\R}^4$. Consider the integral map 
\begin{equation}
F:{\R}^4 \rightarrow {\R}^2: z \mapsto \big( h_1(z), h_2(z) \big) = (c_1, c_2), 
\label{eq-s1one}
\end{equation} 
where $F(0) = (0,0)$.  Suppose that $F$ has the following properties. \medskip 

\indent 1. \parbox[t]{4.25in}{There is an open neighborhood $U$ of the origin $(0,0)$ in ${\R}^2$ such that 
$(0,0)$ is the only critical value of the integral map $F$ in $U$.} 
\smallspace
\indent 2. \parbox[t]{4.25in}{For every $c \in U^{\times } = U \setminus \{ (0,0) \}$ the fiber $F^{-1}(c)$ is noncompact and connected. The fibration 
$\rho = F|F^{-1}(U^{\times }): F^{-1}(U^{\times }) \rightarrow U^{\times }$ is locally trivial.}
\smallspace
\indent 3. \parbox[t]{4.25in}{The singular fiber $F^{-1}(0,0)$ is noncompact and connected. 
For every $z \in F^{-1}(0,0) \setminus \{ 0 \}$ the rank of $DF(z)$ is $2$.} \medskip 

\noindent Then the fibration $\widehat{\rho }= F|F^{-1}(C): F^{-1}(C) \rightarrow C $ 
over the smooth circle $C $ in $U^{\times }$ is \emph{trivial}. So $F^{-1}(C) = 
C \times (S^1 \times \R )$. There is a connection $1$-form $\theta $ on $F^{-1}(C)$, which is 
invariant under the periodic flow of a vector field $X_I$ on $F^{-1}(C)$. For each 
$c  \in C \subseteq U^{\times }$, the curve $t \mapsto {\Gamma }_{\sigma (c)}(t) = 
{\varphi }^{h_1}_t\big( \sigma (c) \big) $, 
where $\sigma $ is a global section of the bundle $\widehat{\rho}$, has scattering phase 
$\Theta (c) = \int_{{\Gamma }_{\sigma (c)}} \theta $. The degree of the map $C \subseteq U^{\times } \rightarrow S^1: c  \mapsto \Theta (c) $ is $-1$. This is the \emph{scattering monodromy} of the focus-focus system. 

\section{The singular fiber}

In this section we show that the singular fiber $F^{-1}(0,0)$ of the integral mapping $F$ (\ref{eq-s1one}) 
is homeomorphic to a once pinched cylinder. \medskip 

Let ${\varphi}^{h_1}_t$ and ${\varphi }^{h_2}_s$ be the flows of the vector 
fields $X_{h_1}$ and $X_{h_2}$, respectively. The hyperbolicity of $X_{h_1}$ at 
$0$ implies that there is an open ball $B$ in ${\R}^4$ (with the Euclidean inner product) centered at $0$ having radius $r$ such that the local stable $W^B_s(0)$ and unstable $W^B_u(0)$ manifolds 
of $0$ in $B$ are smooth connected manifolds, whose tangent space at $0$ is 
the $\mp 1$ eigenspace of the  linear mapping $X_{q_1}$, respectively. The global 
stable manifold $W_s(0)$ of $0$ is ${\bigcup }_{t >0} 
{\varphi }^{h_1}_{-t}(W^B_s(0))$; while the global unstable manifold $W_u(0)$ of $0$ is 
${\bigcup }_{t >0} {\varphi }^{h_1}_t(W^B_u(0))$. When $z \in W_{s,u}(0)$ as $t \rightarrow \infty , \, -\infty $ 
we have $F(z) = F({\varphi }^{h_1}_t(z)) \rightarrow F(0) = (0,0)$. Thus $W_{s,u}(0) \subseteq F^{-1}(0, 0)$. \medskip 

\noindent \textbf{Claim 1.1} $F^{-1}(0,0)\setminus \{ 0 \}  = (W_s(0)\setminus  \{ 0 \} ) \coprod 
(W_u(0) \setminus \{ 0 \}) $. \medskip 

\noindent {\bf Proof}: Since $F^{-1}(0,0)$ is locally invariant under the flow ${\varphi }^{h_1}_t$, it is globally invariant. Thus ${\varphi }^{h_1}_t|F^{-1}(0,0)$ is defined for every $t \in {\R}$. Because of hypothesis $3$, the set $F^{-1}(0,0)^{\times } = F^{-1}(0,0) \setminus \{ 0 \} $ is a smooth $2$-dimensional submanifold of 
${\R}^4$. So ${\varphi }^{h_1}_t(W^B_{s,u}(0)) \setminus \{ 0 \})$ is an open subset of $F^{-1}(0,0)^{\times }$. 
Thus $W_{s,u}(0)^{\times } = W_{s,u}(0) \setminus \{ 0 \} = {\bigcup}_{\mp t \ge 0} 
{\varphi }^{h_1}_t(W^B_{s,u}(0) \setminus \{ 0 \})$. The set $F^{-1}(0,0)$ is invariant 
under the flow ${\varphi }^{h_2}_s$. Because $\{ h_1, h_2 \} =0$, the flows 
${\varphi }^{h_1}_t$ and ${\varphi }^{h_2}_s$ commute. Thus $F^{-1}(0,0)$  is invariant under the 
${\R}^2$-action  
\begin{equation}
\Xi : {\R}^2 \times {\R}^4 \rightarrow {\R}^4: \big( (t,s), z \big) 
\mapsto ({\varphi }^{h_1}_t \comp {\varphi }^{h_2}_s)(z).
\label{eq-s1two}
\end{equation}
So the ${\R}^2$-action ${\Psi }_{(t,s)} = {\Xi}_{(t,s)}|F^{-1}(0,0)$ on $F^{-1}(0,0)$ is defined. 
Because $(0,0) \in {\R}^2$ is an isolated critical value of $F$ by hypothesis $1$, it follows 
that $0 \in {\R}^4$ is an isolated equilibrium point of $X_{h_1}$ and 
$X_{h_2}$. Thus $0$ is an isolated fixed point of the ${\R}^2$-action ${\Psi}_{(t,s)}$ on $F^{-1}(0,0)$. 
If $z \in W_{s,u}(0)$, then ${\varphi }^{h_1}_t \big({\varphi }^{h_2}_s(z) \big) = 
{\varphi }^{h_2}_s \big({\varphi }^{h_1}_t(z) \big) \rightarrow {\varphi }^{h_2}_s(0) =0$ when  
$t \rightarrow \infty , \, -\infty$. So $W_{s,u}(0)$ is invariant under the flow ${\varphi }^{h_2}_s$. 
Because $0$ is a fixed point of the ${\R}^2$-action ${\Psi }_{(t,s)}$ on $F^{-1}(0,0)$, it follows that 
$W_{s,u}(0)^{\times } = W_{s,u}(0) \setminus \{ 0 \} $ is invariant under both flows 
${\varphi }^{h_1}_t$ and ${\varphi }^{h_2}_s$. By hypothesis 2 the vector fields 
$X_{h_1}$ and $X_{h_2}$ are linearly independent 
at each point of $F^{-1}(0,0)^{\times }$. Consequently, every orbit $\mathcal{O}$ of the 
${\R}^2$-action ${\Psi }_{(t,s)}$ on $F^{-1}(0,0)^{\times}$ is open. Because the complement 
of $\mathcal{O}$ in a connected component of $F^{-1}(0,0)^{\times }$ is the union of other ${\R}^2$ orbits of 
${\Psi}_{(t,s)}$, it is also open. Thus $\mathcal{O}$ is a connected component of $F^{-1}(0,0)^{\times }$. 
A similar argument shows that $W_{s,u}(0)^{\times }$ is an ${\R}^2$-orbit in $F^{-1}(0,0)^{\times }$. The orbit 
$\mathcal{O}$ is open and closed in $F^{-1}(0,0)^{\times }$, which implies that it is 
open in $F^{-1}(0,0) = F^{-1}(0,0)^{\times } \cup \{ 0 \} $. If $0$ is not in the closure of 
$\mathcal{O}$ in ${\R}^4$, then $\mathcal{O}$ is closed in $F^{-1}(0,0)$. Hence 
$\mathcal{O}$ is a connected component of $F^{-1}(0,0)$. But $0 \in F^{-1}(0,0)$, 
which is a contradiction. So $\mathcal{O} \cup \{ 0 \}$ is a closed subset of ${\R}^4$. \medskip 

Thus $0$ is the unique limit point in ${\R}^4 \setminus \mathcal{O}$ of the ${\R}^2$-orbit 
$\mathcal{O}$. For any $\mathcal{O}$ in $F^{-1}(0,0)^{\times }$, we know that 
the closure of $\mathcal{O}$ in ${\R}^4$ is $\mathcal{O} \cup \{ 0 \}$. In particular, 
this holds when $\mathcal{O} = W_{s,u}(0)^{\times}$. If $\mathcal{O}$ is an orbit of the 
${\R}^2$-action ${\Psi}_{(t,s)}$ on $F^{-1}(0,0)$ and $z \in \mathcal{O}$, then the mapping 
$(t,s) \mapsto {\varphi }^{h_1}_t({\varphi }^{h_2}_s(z))$ induces a diffeomorphism of 
${\R}^2/J_z$ onto $\mathcal{O}$, where $J_z = \{ (t,s) \in {\R}^2 \, \, \setrule \, \, 
{\varphi }^{h_1}_t({\varphi }^{h_2}_s(z)) = z \} $ is the isotropy group of $z$. $J_z$ is an 
additive subgroup of ${\R}^2$, which does not depend on $z$, because $\mathcal{O}$ is connected. 
Therefore we will write $J_{\mathcal{O}}$ instead of $J_z$. Suppose that 
$\mathcal{O}$ is an ${\R}^2$-orbit of the action ${\Psi}_{(t,s)}$ on $F^{-1}(0,0)^{\times }$ and 
that $J_{\mathcal{O}} \cap ( {\R} \times \{ 0 \} ) \ne \varnothing $. 
Then the flow ${\varphi }^{h_1}_t$ of $X_{h_1}$ would be periodic with period $t_c >0$. 
Because periodic integral curves of $X_{h_1}$, which lie in $\mathcal{O}$ and start 
near $0$ leave a fixed neighborhood of $0$, have an arbitarily large period, we 
deduce that a periodic solution of $X_{h_1}$, which starts near $0$, must stay close to 
$0$. Because $X_{h_1}$ is hyperbolic at $0$ it does not have any periodic 
solutions which remain close to $0$ other than $0$. This is a contradiction. So 
$J_{\mathcal{O}} \cap ( {\R} \times  \{ 0 \}  ) = \varnothing $. \medskip

Combined with the fact that $0$ is the only limit point in ${\R}^4 \setminus \mathcal{O}$ of 
$\mathcal{O}$ and that $\mathcal{O}$ is contained in a connected component of $F^{-1}(0,0)$ of 
${\R}^4$, it follows that for every $z \in \mathcal{O}$ we have 
${\varphi }^{h_1}_t(z) \rightarrow 0$ as $t \rightarrow \infty $ or as $t \rightarrow -\infty $. 
In other words, $z \in {W_s(0)}^{\times }$ or ${W_u(0)}^{\times }$. Because $W_{s,u}(0)^{\times }$ are 
${\R}^2$-orbits of the action ${\Psi}_{(t,s)}$ on $F^{-1}(0,0)^{\times }$, 
it follows that  $F^{-1}(0,0)^{\times } = W_s(0)^{\times } \coprod W_u(0)^{\times }$. So 
$F^{-1}(0,0) = W_s(0) \cup W_u(0)$ and the connected components of 
${F^{-1}(0,0)}^{\times }$ are ${W_{s,u}(0)}^{\times}$. \hfill $\square $ \medskip 
 
\noindent We now prove \medskip 

\noindent \textbf{Claim 1.2 } In a suitable open neighborhood of $0$ in 
$\big( {\bf R}^4, {\omega } =-\dee \alpha = - \dee \, ( p_x \dee x + p_y \dee y) \big)$ there is a 
Hamiltonian function $I$, which equals $h_2 +\mathcal{O}(2)$, whose associated Hamiltonian 
vector field $X_I$ has a periodic flow ${\varphi}^I_u$ and Poisson commutes with the Hamiltonians 
$h_1$ and $h_2$. \medskip 

To construct the function $I$ we prove \medskip 

\noindent \textbf{Lemma 1.3} There is a local diffeomorphism $\Phi $ of 
${\R }^4$, which fixes the origin, is near the identity, and is isotopic to the 
identity map, such that ${\Phi }^{\ast }h_i = q_i $ for $i=1,2$. \medskip %

\noindent \textbf{Proof.} We will use a complex version of the Morse lemma, which is 
proved in the appendix, to construct the desired local diffeomorphism. Introduce complex 
coordinates $(z_1,z_2) = (x- \mathrm{i}y, p_x + \mathrm{i}p_y)$ on ${\R }^4$. Then 
\begin{equation}
\iota : {\R }^4 \rightarrow {\C}^2: (x,y,p_x,p_y) \mapsto (z_1, z_2)
\label{eq-s1twostar}
\end{equation}
is an invertible real linear mapping. The integral map 
\begin{displaymath}
\begin{array}{l}
F:{\R }^4 \rightarrow {\R }^2: \\
\rule{0pt}{12pt}\hspace{.15in} (x,y,p_x,p_y) \mapsto 
\big( xp_x +yp_y +r_1(x,y,p_x,p_y), xp_y-yp_x + r_2(x,y,p_x,p_y) \big) 
\end{array}
\end{displaymath}
becomes the differentiable function
\begin{displaymath}
\mathcal{H}: {\C}^2 \rightarrow \C : (z_1,z_2) \mapsto z_1z_2 + R(z_1,z_2), 
\end{displaymath}
where $R(z_1,z_2) = (r_1 + \mathrm{i} r_2) \big( {\iota }^{-1}(z_1,z_2) \big) $. Because 
$r_i$ is flat to second order at $0$ for $i=1,2$, the function $R$ is flat to second order 
at $(0,0)$. We have $\mathcal{H} = j \comp F \comp {\iota }^{-1}$, where $j: {\R }^2 \rightarrow \C: 
(x,y) \mapsto x + \mathrm{i}y$. Check: $(F \comp {\iota }^{-1})(z_1,z_2) = F(x,y,p_x,p_y)$. So 
\begin{align*}
j\big( (F \comp {\iota }^{-1})(z_1,z_2) \big) & = 
\big( (xp_x+yp_y) + \mathrm{i}(xp_y-yp_x) + (r_1 + \mathrm{i}r_2)(x,y,p_x,p_y) \big) \\
& = z_1z_2 + R(z_1,z_2) .
\end{align*}
Since $D^2\mathcal{H}(0,0) =${\tiny \raisebox{1pt}{$\begin{pmatrix} 0 & 1 \\ 1 & 0 \end{pmatrix}$}} is 
invertible, $(0,0)$ is a nondegenerate critical point of $\mathcal{H}$. By the complex Morse lemma, 
there is an open neighborhood $U$ of $(0,0)$ in ${\C}^2$ and a complex diffeomorphism 
$\varphi : U \rightarrow U$, which fixes $(0,0)$, is near ${\mathrm{id}}_U$, and is isotopic to 
${\mathrm{id}}_U$, such that for every $(z_1, z_2) \in U$
\begin{displaymath}
(\mathcal{H} \comp \varphi ) (z_1,z_2) = \onehalf D^2\mathcal{H}(0,0)\big( (z_1, z_2), (z_1,z_2) \big) 
= \widetilde{\mathcal{H}}(z_1,z_2), 
\end{displaymath}
where $\widetilde{\mathcal{H}}: U \subseteq {\C}^2 \rightarrow \C: (z_1,z_2) \mapsto z_1z_2$. 
In real terms $\widetilde{\mathcal{H}}$ is the integral map $\widetilde{F} = j^{-1} \comp 
\widetilde{\mathcal{H}} \comp \iota $, where
\begin{displaymath}
\widetilde{F}: \widetilde{U} = {\iota }^{-1}(U) \subseteq {\R}^4 \rightarrow {\R }^2: 
(x,y,p_x,p_y) \mapsto ( xp_x + yp_y, xp_y-yp_x) = (q_1,q_2), 
\end{displaymath}
is the integral map of the standard focus-focus system $(q_1, q_2, {\R }^4, \omega )$. 
Also in real terms, the complex 
local diffeomorphism $\varphi $ corresponds to the real local diffeomorphism 
$\Phi = {\iota }^{-1} \comp \varphi \comp \iota $ of $\widetilde{U}$ into itself, which fixes $0$, 
is near the identity, that is, $\Phi = {\mathrm{id}}_U + \mathcal{O}(2)$, and is isotopic to the identity. 
Since $\widetilde{F} = j^{-1} \comp ( \mathcal{H} \comp \varphi ) \comp \iota = 
F \comp \Phi $, we obtain ${\Phi }^{\ast }h_i = q_i$ for $i=1,2$. Warning: $\Phi $ is \emph{not} 
a symplectic diffeomorphism of $(\widetilde{U}, \omega | \widetilde{U})$ into itself. \hfill $\square $ \medskip 

We now begin the construction of the function $I$ in claim 1.2. Let $B$ be 
a ball of radius $r$ in ${\R }^4$ centered at $0$, which is contained in the open set $\widetilde{U}$. Let $Y$ be the vector field ${\Phi }^{\ast }X_{q_2}$ whose flow is ${\psi }_s = {\Phi }^{-1} \comp {\varphi }^{q_2}_s \comp \Phi $. Hence the integral curve ${\Gamma }_w$ of $Y$ starting at $w = {\Phi }^{-1}(z) \in B\setminus \{ 0 \}$ is periodic of period $2\pi $, because 
\begin{displaymath}
{\Gamma }_w(s) = {\psi }_s(w) = {\Phi }^{-1}\big( {\varphi }^{q_2}_s(\Phi (w)) \big) = 
{\Phi }^{-1}({\gamma }_z(s)), 
\end{displaymath}
where ${\gamma }_z$ is an integral curve of $X_{q_2} = 
-y \frac{\partial }{\partial x} +x \frac{\partial }{\partial y}-p_y \frac{\partial }{\partial p_x}+
p_x \frac{\partial }{\partial p_y}$ starting at $z \ne 0$, which is periodic of period $2\pi $. Since 
\begin{displaymath}
L_{Y}h_i = L_{{\Phi }^{\ast }X_{q_2}} h_i = {\Phi }^{\ast }\big( L_{X_{q_2}} ({\Phi }^{-1})^{\ast }h_i \big) 
= {\Phi }^{\ast }(L_{X_{q_2}} q_i ) =0, 
\end{displaymath}
the flow ${\psi }_s$ of $Y$ preserves the level sets of the integral map $F$ 
(\ref{eq-s1one}).\medskip 

For $w \in B\setminus \{ 0 \}$ let  
\begin{equation}
I:B \setminus \{ 0 \} \subseteq {\R}^4 \rightarrow {\R}: w \mapsto I(w) = 
\frac{1}{2\pi} \int_{{\Gamma }_w} \alpha . 
\label{eq-s1three}
\end{equation}  
Then $I = {\Phi }^{\ast }K$, where 
\begin{equation}
K: {\Phi}(B) \setminus \{ 0 \} \subseteq {\R }^4 \rightarrow \R : 
z \mapsto \frac{1}{2\pi} \int_{{\gamma }_z}({\Phi }^{-1})^{\ast } \alpha .
\label{eq-s1four}
\end{equation} 

\noindent {\bf Proof.} We compute 
\begin{align}
I(w) & = \frac{1}{2\pi} \int_{{\Gamma }_w} \alpha = \frac{1}{2\pi} \int^{2\pi}_0 \langle \alpha 
\, \, \setrule \, \, \frac{\dee {\psi }_t}{\dee t} \rangle ({\Gamma }_w(t)) \, \dee t \notag \\
& = \frac{1}{2\pi} \int^{2\pi }_0 \langle \alpha ({\psi }_t(w)) \, \, \setrule \, \, Y({\psi }_t(w)) \rangle  
\, \dee t \notag \\
& = \frac{1}{2\pi} \int^{2\pi }_0 \langle \alpha \big( {\Phi }^{-1}({\varphi }^{q_2}_t(z)) \big) \, \, \setrule 
\, \, T\Phi \, X_{q_2}({\varphi }^{q_2}_t(z)) \rangle \, \dee t, \, \, \, \mbox{since $Y = {\Phi }^{\ast }X_{q_2}$} 
\notag \\
& = \frac{1}{2\pi} \int^{2\pi }_0 \langle ({\Phi }^{-1})^{\ast }\alpha \, \, \setrule \, \, 
X_{q_2} \rangle ({\varphi }^{q_2}_t(z)) \rangle \, \dee t  \notag \\
& = \frac{1}{2\pi} \int^{2\pi }_0 \langle ({\Phi }^{-1})^{\ast }\alpha \, \, \setrule \, \, 
\frac{\dee {\varphi }^{q_2}_t}{\dee t}) \rangle (z) \, \dee t \notag \\
& = \frac{1}{2\pi} \int_{{\gamma }_z}  ({\Phi }^{-1})^{\ast }\alpha = K(z) = K(\Phi (w)). \tag*{$\square $} 
\end{align} 

Next we show that $I$ (\ref{eq-s1three}) is smooth near $0$. \medskip 
 
\noindent {\bf Proof.} Let $z_1 = x-i y$ and $z_2 = p_x+i p_y$. 
Then $q_2 = \mathrm{Im}\, z_1z_2$ and the flow ${\varphi }^{q_2}_t$ of $X_{q_2}$ is 
$\big( t, (z_1, z_2) \big) \mapsto ({\mathrm{e}}^{it}z_1, {\mathrm{e}}^{-it}z_2)$. Let $D = 
\{ \zeta \in \C \, \, \setrule \, \, |\zeta | \le 1 \} $ with boundary 
$\partial D = S^1 = \{ \zeta \in \C \, \, \setrule \, \, |\zeta | = 1 \}$. Let 
\begin{displaymath}
k: D \times {\C}^2 \rightarrow {\C}^2: \big( \zeta , (z_1,z_2) \big) 
\mapsto (\zeta z_1, {\zeta }^{-1} z_2 )
\end{displaymath}
with $k_z: D \rightarrow {C}^2: \zeta \mapsto (\zeta z_1, {\zeta }^{-1} z_2)$. 
Using Stokes' theorem we have 
\begin{displaymath}
K(z) = \int_{\partial D} k^{\ast}_z ({\Phi }^{-1})^{\ast }\alpha = 
-\int_D k^{\ast }_z({\Phi }^{-1}(\omega )) . 
\end{displaymath}
Since $D$ is compact and $\omega $, $k_z$ are smooth, 
it follows that $K$ is smooth near $0$. Thus $I = {\Phi }^{\ast }K$ is smooth near $0$. \hfill $\square $ \medskip
 
\noindent \textbf{Claim 1.4} The function $I$ (\ref{eq-s1three}) is an action for 
$\big( h_1, h_2, {\R}^4, \omega \big)$. \medskip 

This follows from the next three results. \medskip 

\noindent 1. The function $I$ Poisson commutes with $h_i$ for $i=1,2$ on $F^{-1}(c)$, 
where $c$ is a regular value of the integral map $F$ (\ref{eq-s1one}). \medskip 

\noindent {\bf Proof.} We compute. $\{ I, h_i \} = L_{X_{h_i}}I = \int_{{\Gamma }_w} L_{X_{h_i}} \alpha $, 
because we can move ${\Gamma }_w$ by a homotopy in 
$F^{-1}(c)$ without changing the integral, the new integral does not 
depend on $w$. Thus we can take the Lie derivative under the integral sign. But
\begin{align} 
\int_{{\Gamma }_w} L_{X_{h_i}} \alpha  & = 
\int_{{\Gamma }_w} X_{h_i} \, \lefthook \, \dee \alpha + \dee \, (X_{h_i} \, \lefthook \, \alpha ) 
\notag \\
&\hspace{-.5in}= \int_{{\Gamma }_w} \dee \, (- h_i + X_{h_i} \, \lefthook \, \alpha ) =0, \, \, 
\parbox[t]{2in}{since ${\Gamma }_w$ is a closed curve.} \tag*{$\square $}
\end{align}

\noindent 2. For all values of $c$ close to but not equal to $0$ and for all $w \in F^{-1}(c)$, the tangent vectors $X_I(w)$ and $X_{h_1}(w)$ to $F^{-1}(c)$ at $w$ are linearly independent. \medskip

\vspace{-.15in}\noindent {\bf Proof.} Since $\Phi = {\mathrm{id}}_U + 
{\mathcal{O}}(1)^2$, we get ${\Phi }^{-1} = {\mathrm{id}}_U + {\mathcal{O}}(1)^2$. So 
$({\Phi }^{-1})^{\ast } \alpha = \alpha + \mathcal{O}(1)$, which gives 
\begin{equation}
K(z) = \frac{1}{2\pi } \, \int_{{\gamma }_z} ({\Phi }^{-1})^{\ast } \alpha = 
\Big( \frac{1}{2\pi } \, \int_{{\gamma }_z} \alpha \Big) + \mathcal{O}(2) = q_2 + \mathcal{O}(2). 
\label{eq-s1five}
\end{equation}
Therefore 
\begin{equation}
I = {\Phi }^{\ast }K = {\Phi }^{\ast }q_2 + \mathcal{O}(2) = h_2 + \mathcal{O}(2). 
\label{eq-s1six}
\end{equation}
Since the vector fields $X_{h_1}$ and $X_{h_2}$ are linearly independent on $F^{-1}(c)$ for 
$c$ near, but not at the origin, the vector fields $X_I$ and $X_{h_1}$ are also.  \hfill $\square $ \medskip 

\noindent 3. For all $c$ close to but not equal to $0$, the flow ${\varphi}^I_u$ of $X_I$ on $F^{-1}(c)$ is periodic of period $T_c$. \medskip 

\noindent {\bf Proof.} For all $c \in {\R }^2$ near but not at $0$, there are smooth functions 
$a(c)$ and $b(c)$ such that 
\begin{displaymath}
X_I = a(c) X_{h_1} + b(c)X_{h_2}
\end{displaymath}
on $F^{-1}(c)$. Thus ${\varphi }^I_t ={\varphi }^{h_1}_{a(c)t} \comp {\varphi }^{h_2}_{b(c)t}$ is 
defined for all $t \in \R$, since the vector fields $X_{h_i}$ for $i=1,2$ are complete. Hence 
the vector field $X_I$ on $F^{-1}(c)$ is complete. Since ${\Phi }_{\ast }I = K$, 
the vector field ${\widetilde{X}}_K$ on $\big( {\R }^4, 
\widetilde{\omega } = {\Phi }_{\ast}\omega \big)$, which equals ${\Phi }_{\ast }X_I$, is complete. \medskip 

From $0 = \{ I, h_i \} = \{ {\Phi }^{\ast }K, {\Phi }^{\ast }q_i \} = {\Phi}^{\ast } \{ \{ K, q_i \} \} $, 
where $\{ \{ \, \, , \, \, \} \}$ is the Poisson bracket on $({\R }^4, \widetilde{\omega })$, it follows that 
the vector fields ${\widetilde{X}}_K$ and ${\widetilde{X}}_{q_i}$ on $({\R }^4, \widetilde{\omega })$ 
commute. So their flows ${\widetilde{\varphi }}^K_u$ and ${\widetilde{\varphi }}^{\, q_i}_t$ commute. 
Recall that $\widetilde{F}: {\R }^4 \rightarrow{\R }^2: z \mapsto \big( q_1(z), q_2(z) \big) $. 
The flow ${\widetilde{\varphi }}^K_u$ on ${\widetilde{F}}^{-1}(c)$ is periodic for every 
$c \in {\R}^2$ close to but not at $0$. To see this we argue as follows. \medskip 

We have an ${\R }^2$ action 
\begin{displaymath}
\Lambda : {\R }^2 \times {\widetilde{F}}^{-1}(c) \rightarrow {\widetilde{F}}^{-1}(c): 
\big( (u,t), z \big) \mapsto ({\widetilde{\varphi}}^K_u \comp {\widetilde{\varphi }}^{\, q_1}_t)(z).  
\end{displaymath}
For any $\widetilde{w} \in {\widetilde{F}}^{-1}(c)$ the isotropy group 
${\Lambda}_{\widetilde{w}} = \{ (u,t) \in {\R }^2 \setrule \, \Lambda (u,t, \widetilde{w}) = \widetilde{w} \} $ 
is a rank $1$ lattice, since ${\R }^2/{\Lambda}_{\widetilde{w}} = {\widetilde{F}}^{-1}(c) = S^1 \times \R $. 
Let $(u_0,t_0) \in {\Lambda}_{\widetilde{w}}$. Then $({\widetilde{\varphi}}^K_{u_0} \comp 
{\widetilde{\varphi }}^{\, q_1}_{t_0}) (\widetilde{w}) = \widetilde{w}$. Suppose that there is $t' >0$ such that 
$(u_0, t_0+t') \in {\Lambda}_{\widetilde{w}}$. Then $({\widetilde{\varphi}}^K_{u_0} \comp 
{\widetilde{\varphi }}^{\, q_1}_{t_0+t'}) (\widetilde{w}) = \widetilde{w}$. So 
${\widetilde{\varphi }}^{\, q_1}_{-t'}(\widetilde{w}) = {\widetilde{\varphi}}^K_{u_0} 
\big( {\widetilde{\varphi }}^{\, q_1}_{t_0} (\widetilde{w}) \big) $, since ${\widetilde{\varphi}}^K_{u}$ and 
${\widetilde{\varphi }}^{\, q_1}_{t}$ commute. Thus ${\widetilde{\varphi }}^{\, q_1}_{t} \big( 
{\widetilde{\varphi }}^{\, q_1}_{-t'}(\widetilde{w}) \big) = {\widetilde{\varphi }}^{\, q_1}_{t}(\widetilde{w})$, 
which implies that the integral curve $t \mapsto  {\widetilde{\varphi }}^{\, q_1}_{t}(\widetilde{w})$ of 
${\widetilde{X}}_{q_1}$ is periodic of period $-t'$. Since the diffeomorphism $\Phi $ is isotopic 
to the identity map, each integral curve of ${\widetilde{X}}_{q_1}$ is homotopic to an 
integral curve of $X_{q_1}$, because the symplectic form $\widetilde{\omega }$ is homotopic to the 
symplectic form $\omega $. See lemma 4.3. But $X_{q_1}$ has no periodic integral curves. 
Thus $t' =0$. Since the rank of ${\Lambda }_{\widetilde{w}}$ is $1$, 
there is a $u' >0$ such that $(u_0 +u',t_0) \in 
{\Lambda }_{\widetilde{w}}$, that is, $({\widetilde{\varphi }}^K_{u'} \comp {\widetilde{\varphi }}^K_{u_0}) 
\big( {\widetilde{\varphi }}^{q_1}_{t_0}(\widetilde{w}) \big) = \widetilde{w}$. So 
${\widetilde{\varphi }}_{u'}(\widetilde{w}) = \widetilde{w}$, 
since $({\widetilde{\varphi }}^K_{u_0} \comp {\widetilde{\varphi }}^{\, q_1}_{t_0})(\widetilde{w}) = 
\widetilde{w}$. Thus the integral curve $u \mapsto {\widetilde{\varphi }}^K_u(w)$ of ${\widetilde{X}}_K$ is 
periodic of period $u'$. Since $\widetilde{w}$ is an arbitrary point of ${\widetilde{F}}^{-1}(c)$, the 
flow ${\widetilde{\varphi }}^K_u$ of ${\widetilde{X}}_K$ on ${\widetilde{F}}^{-1}(c)$ is periodic of 
period $T_c = u'$. Thus the flow of $X_I$ on $F^{-1}(c)$ is periodic of period $T_c$, because 
$I = {\Phi }^{\ast }K$. This proves 3 and completes the proof of claim 1.4. \hfill $\square $ \medskip 

This completes the proof of claim 1.2. \hfill $\square $ \medskip  

\noindent \textbf{Claim 1.5} $F^{-1}(0,0)$ is homeomorphic to a pinched cylinder, that is, 
a cylinder $S^1 \times \R $ with one of its generating circles pinched to the origin $0$. The 
singular fiber $F^{-1}(0,0)$ has two transverse tangent planes at $0$. \medskip 

\noindent {\bf Proof.} Since the action $I$ Poisson 
commutes with the integrals $h_i$ for $i=1,2$, the flow ${\varphi }^I_u$ of $X_I$ leaves the 
fiber $F^{-1}(0,0) \cap V$ invariant. Here $V \subseteq B$ is an open neighborhood of $0$, 
which is invariant under the $S^1$-action generated by ${\varphi}^I_u$. 
Note that $W_{s,u}(0) \cap V$ is invariant under the flow ${\varphi }^I_u$. \medskip 

We now extend the $S^1$-action ${\varphi }^I_u|(W_{s,u}(0) \cap V)$ to all of $W_{s,u}(0)$. Let 
$p \in W_s(0)$ or $W_u(0)$. Then there is an open neighborhood $V_p$ of $0$ in ${\bf R}^4$ and a 
time $t_p > 0$ such that ${\varphi }^{h_1}_{\mp t_p}(V_p) \subseteq V$. For every $\widetilde{p} 
\in V_p$ let ${\widehat{\varphi }}^I_u(\widetilde{p}) = ({\varphi }^{h_1}_{\pm t_p} \comp {\varphi }^I_u 
\comp {\varphi }^{h_1}_{\mp t_p})(\widetilde{p})$, where {\tiny $\left\{ \begin{array}{rl} -, & \mbox{if
$p \in W_s(0)$} \\ +, & \mbox{if $p \in W_u(0)$.} \end{array} \right. $} Then ${\widehat{\varphi }}^I_u$ defines an $S^1$-action on an open neighborhood $\mathcal{V} = 
\bigcup_{p \in W_s(0)} V_p \cup \bigcup_{p \in W_u(0)} V_p$ of $F^{-1}(0,0) = W_s(0) \cup W_u(0)$. 
Therefore we have an ${\R}^2$-action on $\mathcal{V}$ defined by 
\begin{equation}
\widetilde{\Xi }: {\R}^2 \times \mathcal{V} \rightarrow \mathcal{V}: 
\big( (u,t), z \big) \mapsto {\widehat{\varphi }}^I_u \comp {\varphi }^{h_1}_t(z). 
\label{eq-s1seven}
\end{equation}
Note that $F^{-1}(0,0)$ is invariant under the action $\widetilde{\Xi}$ and that $0$ is 
the only fixed point of the flow ${\varphi }^I_u$ on $F^{-1}(0,0)$. So we 
have an ${\bf R}^2$-action ${\widetilde{\Psi}}_{(u,t)}= \widetilde{\Xi}|F^{-1}(0,0)^{\times }$. 
Let $J_{\mathcal{O}}$ be the isotropy group for the ${\bf R}^2$-orbit 
$\mathcal{O} = {W_s(0)}^{\times}$ or ${W_u(0)}^{\times }$. Because $F^{-1}(0)^{\times }$ is not compact, the rank of the lattice $J_{\mathcal{O}}$ can not be equal to $2$. But $J_{\mathcal{O}} \ne 
\varnothing $. So $J_{\mathcal{O}}$ is isomorphic to $\Z$. 
Since the flow ${\varphi }^I_u|F^{-1}(0,0)$ is periodic of period $T_{(0,0)} >0$, it follows that 
$J_{\mathcal{O}} = T_{(0,0)} \Z$. Thus $\mathcal{O}$ is diffeomorphic 
to the cylinder ${\R}^2/J_{\mathcal{O}} = (\R /T_{(0,0)} \Z) \times \R = 
S^1 \times {\R}$. Because $F^{-1}(0,0) = F^{-1}(0,0)^{\times } \cup \{ 0 \} $ and $\{ 0 \} $ is the 
only limit point of the closure of $F^{-1}(0,0)^{\times } = \mathcal{O}$ in ${\R}^4$, it follows that 
$F^{-1}(0,0)$ is homeomorphic to the one point compactification of the cylinders $S^1 \times \R = 
{W_{s,u}(0)}^{\times }$. In other words, $F^{-1}(0,0)$ is 
homeomorphic to a smooth cylinder $S^1 \times \R $ with 
a generating circle pinched to the point $0$. Since the stable and unstable manifolds of the linear vector field $DX_{h_1}(0)$ at $0$ are the coordinate $2$-planes $\{ 0 \} \times {\R}^2$ and ${\R}^2 \times \{ 0 \}$, respectively, in ${\R}^4$, the limits of the tangent planes to $F^{-1}(0,0)^{\times }$ at $0$ exist and are transverse. Thus $F^{-1}(0,0)$ is a smooth $2$-sphere, which is immersed in ${\R}^4$ with 
a normal crossing at $0$. \hfill $\square $  

\section{Nearby regular fibers}

In this section we study the fibers of the integral map $F$ (\ref{eq-s1one}), which are close to the 
singular fiber $F^{-1}(0,0)$. We prove \medskip 

\noindent \textbf{Claim 2.1} There is an open neighborhood $W$ of $F^{-1}(0,0)$ in ${\R }^4$, which is 
invariant under the flows ${\varphi }^{h_1}_t$ and ${\varphi }^I_u$, and an open neighborhood 
$U$ of $(0,0)$ in ${\R }^2$ such that 
\begin{equation}
F|(W \setminus F^{-1}(0,0)): W \setminus F^{-1}(0,0) \rightarrow U \setminus \{ (0,0) \}: 
p \mapsto \big( h_1(p), h_2(p) \big) 
\label{eq-s2one}
\end{equation}
is a smooth surjective submersion, which defines a locally trivial fibration whose fiber $F^{-1}(c) \cap W$ for each $c \in U \setminus \{ (0,0) \}$ is a smooth cylinder $S^1 \times \R $, that is, an orbit of the action 
\begin{displaymath}
\widehat{\Xi}: {\R }^2 \times W \rightarrow W: \big( (u,t), p \big) \mapsto 
({\varphi }^I_u \comp {\varphi }^{h_1}_t)(p). 
\end{displaymath}

\noindent \textbf{Proof.} Let $B$ be an open ball in ${\R }^4$ centered at $0$ whose radius is small enough that $\partial B$ intersects $F^{-1}(0,0)$ in two circles $W^B_{s,u}(0) \cap \partial B$. We can arrange that 
all of the orbits of the ${\R }^2$ action 
\begin{displaymath}
\widetilde{\Xi}: {\R }^2 \times \mathcal{V} \rightarrow \mathcal{V}: 
\big( (u,t),z \big) \mapsto ({\varphi }^I_u \comp {\varphi }^{h_1}_t)(z), 
\end{displaymath}
where $\mathcal{V}$ is an open neighborhood of $F^{-1}(0,0)$ in ${\R }^4$, which is 
invariant under the flow ${\varphi }^I_u$, are $2$ dimensional, and all of its orbits near 
$0$ intersect $\partial B$ in two circles, which are close to the 
circles $W_{s,u}(0) \cap \partial B$. \medskip 

Consider the \emph{local} ${\R }^2$ action ${\widetilde{\Xi}}_{(u,t)}|\overline{B}$. 
Let $W $ be the union of ${\R }^2$-orbits which intersect $\mathcal{V} \cap \partial B$ or 
$\{ 0 \} $. The union of $\widetilde{\Xi}$ orbits in $\overline{B}$ which intersect of $\mathcal{V} \cap \partial B$,  is an open subset of $\overline{B}$, which contains a small 
neighborhood $V$ of $0$ in ${\R }^4$. This follows because by hyperbolicity of the integral curves 
of $X_{h_1}$ at a distance $\delta > 0$ from $0$ enter and leave $B$ at points on $\partial B$, 
which are at a distance $\mathrm{O}(\delta )$ from $W_{s,u}(0) \cap \partial B$. Thus 
$W$ is an open neighborhood of $F^{-1}(0,0)$ in ${\R }^4$ such that $F^{-1}(c) \cap W$ is equal to 
the $\widetilde{\Xi}$ orbit in $\overline{B}$ through $F^{-1}(c) \cap \partial B$ for every $c\in 
U \setminus \{ (0,0) \}$. Hence the smooth mapping 
$F|W: W \rightarrow U\setminus \{ (0,0) \}$ is surjective with connected fibers. The invariance 
of $F$ under the local ${\R }^2$ action ${\widetilde{\Xi}}_{(u,t)}|\overline{B}$, together with the 
fact that at every $z \in \mathcal{V} \cap \partial B$ the rank of $DF(z)$ is $2$, implies that 
at each point $w$ on an ${\R }^2$-orbit which intersects $\mathcal{V} \cap \partial B$, the rank 
of $DF(w)$ is $2$. Thus the map (\ref{eq-s2one}) is a surjective submersion, which defines a fibration with connected fibers. By hypothesis 2 this fibration is locally trivial. \medskip 

We now show that each fiber of the fibration (\ref{eq-s2one}) is a smooth cylinder. Since the 
flows ${\varphi }^{h_1}_t$ and ${\varphi }^I_u$ leave the fibers $F^{-1}(c) \cap W$ invariant, 
they define an ${\R }^2$-action 
\begin{displaymath}
{\Xi}^{\vvee}: {\R }^2 \times (W \setminus F^{-1}(0,0)) \rightarrow W \setminus F^{-1}(0,0): 
\big( (u,t), z \big) \mapsto ({\varphi }^I_u \comp {\varphi }^{h_1}_t)(z)
\end{displaymath}
on $W \setminus F^{-1}(0,0)$. Because the vector fields $X_{h_1}$ and $X_I$ are linearly 
independent at each point of $W \setminus F^{-1}(0,0)$, an ${\R }^2$ orbit $\mathcal{O}$ is 
an open subset of $W \setminus F^{-1}(0,0)$. Since $W \setminus F^{-1}(0,0)$ is connected, 
it follows that $\mathcal{O} =  W \setminus F^{-1}(0,0)$. Now $\mathcal{O}$ is diffeomorphic to 
${\R}^2/J_{\mathcal{O}}$ and $\mathcal{O}$ is not compact. Thus $J_{\mathcal{O}}$ is a 
rank $1$ lattice in ${\R }^2$. Hence $\mathcal{O}$ is a smooth cylinder $S^1 \times \R $. 
\hfill $\square $ \medskip 

Let ${\Sigma }_{\xi}$ be the image of a smooth local section $\sigma : U \setminus \{ (0,0) \} \rightarrow 
W \setminus F^{-1}(0,0)$ of the fibration (\ref{eq-s2one}) at $\xi \in W \setminus F^{-1}(0,0)$, 
which is invariant under the flow ${\varphi }^I_u$ of $X_I$. Because $T_c$ is the period of 
the flow ${\varphi }^I_u$, it is the period of every integral curve $u \mapsto {\varphi }^I_u({\xi }')$ 
for every ${\xi }' \in {\Sigma }_{\xi } \cap F^{-1}(c)$ with $c \in U \setminus \{ (0,0) \}$. Thus for 
every $c \in U \setminus \{ (0,0) \}$ we have $(T_c,0) \in J_{{\xi}'}$. \medskip 

Let $\xi \in W^B_s(0)$ and $\eta \in W^B_u(0)$. Let ${\Sigma }_{\xi , \eta } $ be the image 
in $B$ of a smooth local section of the fibration (\ref{eq-s2one}) at $\xi $, $\eta $, which is 
invariant under the $S^1$ action ${\varphi }^u_I$ on $W$. For each $c \in U \setminus \{ (0,0) \}$, 
it follows that $F^{-1}(c) \cap (W \cap {\Sigma }_{\xi , \eta })$ is diffeomorphic to a circle 
${\mathcal{C}}_{\xi ,\eta }(c)$, which is an orbit of the $S^1$ action ${\varphi }^I_u$. For each 
$\xi (c) \in {\mathcal{C}}_{\xi (c)}$ there is a smallest positive time $\tau (c)$ such that 
the integral curve $t \mapsto {\varphi }^{h_1}_{\tau (c)}\big( \xi (c) \big)$ lies in ${\mathcal{C}}_{\eta }(c)$. 
In other words, $\eta (c) = {\varphi }^{h_1}_{\tau (c)}\big( \xi (c) \big) \in {\mathcal{C}}_{\eta }(c)$. 
As $\xi (c)$ traces out ${\mathcal{C}}_{\xi }(c)$ once, $\eta (c)$ traces out ${\mathcal{C}}_{\eta }(c)$ 
once. Thus the circles ${\mathcal{C}}_{\xi , \eta }(c)$ bound a subset of $F^{-1}(c)$, which is 
diffeomorphic to a compact cylinder $S^1 \times [0,1]$. 

\section{Holomorphic focus-focus system}

In this section we give a complex variables treatment of the standard focus-focus system 
$(q_1, q_2, {\R }^4, \omega )$. \medskip 

Let $(z_1, z_2) = (x-\mathrm{i}y, p_x + \mathrm{i}p_y)$ be coordinates on ${\C}^2$. 
Let $\varpi = \dee z_1 \wedge \dee z_2$ be a complex 
symplectic form on ${\C}^2$. Consider the holomorphic Hamiltonian 
\begin{displaymath}
\mathcal{H}: {\C}^2 \rightarrow \C: (z_1, z_2) \mapsto z_1z_2 = (xp_x +yp_y) + 
\mathrm{i}(xp_y-yp_x)  = q_1 + \mathrm{i}q_2.
\end{displaymath}
The complex Hamiltonian vector field on $({\C}^2, \varpi )$ associated to $\mathcal{H}$ is 
\begin{displaymath}
X_{\mathcal{H}} = \frac{\partial \mathcal{H}}{\partial z_2} \frac{\partial }{\partial z_1} - 
\frac{\partial \mathcal{H}}{\partial z_1} \frac{\partial }{\partial z_2} = 
z_1 \frac{\partial }{\partial z_1} - z_2 \frac{\partial }{\partial z_2}, 
\end{displaymath}
since
\begin{displaymath}
X_{\mathcal{H}} \lefthook (\dee z_1 \wedge \dee z_2) = z_1 \dee z_1 + z_2 \dee z_2 = \dee \, (z_1z_2) 
= \dee \mathcal{H}.
\end{displaymath}
The complex integral curves of $X_{\mathcal{H}}$ satisfy 
\begin{displaymath}
\frac{\dee z_1}{\dee \tau } = z_1 \, \, \, \mathrm{and} \, \, \, \frac{\dee z_2}{\dee \tau } = - z_2. 
\end{displaymath}
Here $\tau $ is \emph{complex} time parameter. The \emph{complex} flow of $X_{\mathcal{H}}$ on ${\C}^2$ is 
\begin{displaymath}
{\varphi }^{\mathcal{H}}: \C \times {\C}^2 \rightarrow {\C}^2: \big( \tau, (z_1, z_2) \big) \mapsto 
({\mathrm{e}}^{\tau}z_1, {\mathrm{e}}^{-\tau }z_2). 
\end{displaymath}

Let 
\begin{displaymath}
\widetilde{\Sigma }: {\C}^{\times } = \C \setminus \{ 0 \} \rightarrow {\C}^2: c \mapsto (c,1). 
\end{displaymath}
Since $\mathcal{H}(\widetilde{\Sigma }(c)) = c$ for every $c \in {\C}^{\times}$, it follows that 
$\widetilde{\Sigma} (c) \in {\mathcal{H}}^{-1}(c)$ for every $c \in {\C}^{\times }$. Thus $\widetilde{\Sigma }$ 
is a global section of the bundle
\begin{displaymath}
\widehat{\rho} = \mathcal{H}|{\mathcal{H}}^{-1}({\C}^{\times }):{\mathcal{H}}^{-1}({\C}^{\times }) 
\rightarrow {\C}^{\times}: (z_1, z_2) \mapsto \mathcal{H}(z_1,z_2). 
\end{displaymath}

\noindent \textbf{Lemma 3.1} The bundle $\widehat{\rho }$ is trivial. \medskip 

\noindent \textbf{Proof.} Consider the map 
\begin{displaymath}
\tau : {\mathcal{H}}^{-1}({\C}^{\times }) \subseteq {\C}^2 \rightarrow {\C}^{\times } \times {\mathcal{H}}^{-1}(1): 
(z_1, z_2) \mapsto \big( z_1z_2, (z_1, z_2(z_1z_2)^{-1}) \big) , 
\end{displaymath}
whose inverse is 
\begin{displaymath}
{\tau }^{-1}: {\C}^{\times } \times {\mathcal{H}}^{-1}(1) \rightarrow {\mathcal{H}}^{-1}({\C}^{\times}): 
\big( c, (w, w^{-1}) \big) \mapsto (w, cw^{-1}). 
\end{displaymath}
Check: 
\begin{align*}
(\tau \comp {\tau }^{-1})\big( c, (w, w^{-1}) \big) & = \tau (w, cw^{-1}) =
\big( c, (w, cw^{-1}(wcw^{-1})^{-1} \big) \\
& = \big( c, (w, w^{-1}) \big) 
\end{align*}
and 
\begin{align*}
({\tau }^{-1} \comp \tau )(z_1, z_2) & = {\tau }^{-1}\big( z_1z_2, (z_1, z_2(z_1z_2)^{-1} ) \big)
= \big( z_1, (z_1z_2) z_2(z_1z_2)^{-1} \big) \\
& = (z_1, z_2). 
\end{align*}
Thus the map $\tau $ trivializes the bundle $\widehat{\rho }$. \hfill $\square $ \medskip 

In real terms, the function $\mathcal{H}:{\C}^2 \rightarrow \C$ is the energy momentum map 
$\widetilde{F}:{\R}^4 \rightarrow {\R }^2$ of the focus-focus system. Here $\widetilde{F} = 
j \comp \mathcal{H} \comp {\iota }^{-1}$, where $\iota $ (\ref{eq-s1twostar}) and 
$j: \C \rightarrow {\R }^2: z \mapsto (\mathrm{Re}\, z, \mathrm{Im}\, z)$. \medskip  

In real terms the section $\widetilde{\Sigma }: {\C}^{\times } \rightarrow {\C}^2:c \mapsto (c,1)$ of 
the bundle $\widehat{\rho }$ is the map $\sigma :({\R}^2)^{\times } = 
{\R }^2\setminus \{ 0 \} \rightarrow {\R}^4$, where 
$\sigma = {\iota }^{-1}\comp \widetilde{\Sigma} \comp j^{-1}$. We calculate $\sigma $. For $(c_1, c_2) \in 
({\R }^2)^{\times }$, we have $j^{-1}(c_1,c_2) = c_1 + \mathrm{i}c_2$. So 
$(\widetilde{\Sigma } \comp j^{-1})(c_1,c_2) = (c_1+\mathrm{i}c_2,1)$, which implies 
${\sigma }(c_1,c_2) = {\iota }^{-1}(c_1 + \mathrm{i}c_2,1) = (c_1, -c_2,1,0)$. \medskip 

\noindent \textbf{Lemma 3.2} The fibration
\begin{displaymath}
\widetilde{F}|{\widetilde{F}}^{-1}(R): {\widetilde{F}}^{-1}(R) \subseteq {\R}^4 \rightarrow R = ({\R}^2)^{\times} \subseteq {\R }^2: z \mapsto \big( q_1(z), q_2(z) \big) 
\end{displaymath}
is trivial. \medskip 

\noindent \textbf{Proof.} Since $\widetilde{F} \big( \sigma (c_1,c_2) \big)  = \widetilde{F}(c_1,-c_2,1,0) = (c_1,c_2)$ for every $(c_1, c_2) \in R$, the map $\sigma $ is a global section of the bundle 
$\widetilde{F}|{\widetilde{F}}^{-1}(R)$. So ${\widetilde{F}}^{-1}(R)$ is diffeomorphic to 
$R \times {\widetilde{F}}^{-1}(1,0)$, where ${\widetilde{F}}^{-1}(1,0)$ is a cylinder $S^1 \times \R $. 
\hfill $\square $ \medskip 

In real terms the complex flow ${\varphi }^{\mathcal{H}}_{\tau }$ of the Hamiltonian vector field 
$X_{\mathcal{H}}$ is ${\iota }^{-1} \comp {\varphi }^{\mathcal{H}}_{\tau } \comp \iota $, where $\tau = t + \mathrm{i}s$. We compute. 
\begin{align*}
({\iota }^{-1} \comp {\varphi }^{\mathcal{H}}_{\tau } \comp \iota )(x,y,p_x,p_y) & = 
{\iota }^{-1}\big( {\mathrm{e}}^{t+\mathrm{i}s}(x-\mathrm{i}y), 
{\mathrm{e}}^{-t-\mathrm{i}s}(p_x+\mathrm{i}p_y) \big)  \\ 
& \hspace{-1in}= {\iota }^{-1}\big( {\mathrm{e}}^t\big[ (x\cos s +y \sin s) +\mathrm{i}(x \sin s - y \cos s) \big] , \\
& \hspace{-.5in}{\mathrm{e}}^{-t}\big[ (p_x\cos s +p_y \sin s) +\mathrm{i}(-p_x\sin s + p_y \cos s) \big] \big) \\
&\hspace{-1in} = \big( {\mathrm{e}}^t(x\cos s + y \sin s), {\mathrm{e}}^t(-x\sin s + y \cos s), \\
&\hspace{-.5in} {\mathrm{e}}^{-t}(p_x\cos s + p_y \sin s), {\mathrm{e}}^{-t}(-p_x\sin s + p_y \cos s) \big) \\
& \hspace{-1in} = {\varphi }^{q_1}_t\Big( \mbox{\footnotesize $\begin{pmatrix} \cos s & \sin s \\ -\sin s & \cos s \end{pmatrix} \, \begin{pmatrix} x \\ y \end{pmatrix}$} , 
\mbox{\footnotesize $\begin{pmatrix} \cos s & \sin s \\ -\sin s & \cos s \end{pmatrix} \,  
\begin{pmatrix} p_x \\ p_y \end{pmatrix}$} \Big) \\
& \hspace{-1in}= ({\varphi }^{q_1}_t \comp {\varphi }^{q_2}_{-s})(x,y,p_x,p_y). 
\end{align*}

Fix $\varepsilon >0$. Consider the complex curves 
\begin{subequations}
\begin{align}
&\xi : \C \rightarrow {\C}^2: c \mapsto (\ttfrac{c}{\varepsilon }, \varepsilon ) 
\label{eq-s3onea} \\
&\hspace{-1.7in}\textrm{and} \notag \\ 
& \eta : \C \rightarrow {\C}^2: c \mapsto (\varepsilon , \ttfrac{c}{\varepsilon}). 
\label{eq-s3oneb}
\end{align}
\end{subequations}

\noindent \textbf{Lemma 3.3} The curves $\xi $ and $\eta $ are transverse to ${\mathcal{H}}^{-1}(c)$ for 
every $c \in {\C}^{\times}$. \medskip 

\vspace{-.15in}\noindent \textbf{Proof.} We have $T_{c}\xi = \frac{1}{\varepsilon} \frac{\partial }{\partial z_1}$;  
while $T_{(z_1,z_2)}{\mathcal{H}}^{-1}(c) = z_1\frac{\partial }{\partial z_1} + z_2\frac{\partial }{\partial z_2}$. 
Since $(z_1, z_2) \in {\mathcal{H}}^{-1}(c)$, both $z_1$ and $z_2$ are nonzero. Thus $T_{c} \xi$ does not 
lie in $T_{(z_1,z_2)}{\mathcal{H}}^{-1}(c)$. Hence the curve $\xi $ (\ref{eq-s3onea}) is transverse to 
${\mathcal{H}}^{-1}(c)$ in ${\C}^2$ for every $c \ne 0$. A similar argument shows that the curve $\eta $ is tranverse to ${\mathcal{H}}^{-1}(c)$ for every $c \ne 0$. \hfill $\square $ \medskip 

Fix $c \in \C$. Look at the circles 
\begin{subequations}
\begin{align}
& {\xi }_s: \R \rightarrow {\mathcal{H}}^{-1}(c): s \mapsto (\ttfrac{c}{{\mathrm{e}}^{-\mathrm{i}s}\varepsilon }, 
{\mathrm{e}}^{-\mathrm{i}s}\varepsilon ) 
\label{eq-s3twoa} \\
&\hspace{-1.3in} \textrm{and} \notag \\
& {\eta }_s: \R \rightarrow {\mathcal{H}}^{-1}(c): s \mapsto ({\mathrm{e}}^{\mathrm{i}s} \varepsilon , 
\ttfrac{c}{{\mathrm{e}}^{\mathrm{i}s}\varepsilon }). 
\label{eq-s3twob}
\end{align}
\end{subequations}

\noindent \textbf{Lemma 3.4} For every $s \in \R$ the point ${\xi }_s(0) = (0, {\mathrm{e}}^{-\mathrm{i}s}\varepsilon) $ lies on the stable manifold $W_s(0)$ of the vector field $X_{q_2}$ on ${\R }^4$; while the point 
${\eta }_s(0) = ( {\mathrm{e}}^{\mathrm{i}s}\varepsilon, 0)$ lies on the unstable manifold $W_u(0)$ 
of the vector field $X_{q_1}$ for every $s \in \R $. 
\medskip 

\noindent \textbf{Proof.} We compute. 
\begin{displaymath}
{\varphi }^{\mathcal{H}}_{t + \mathrm{i}u}(0, {\mathrm{e}}^{-\mathrm{i}s}\varepsilon) = 
\big( 0, {\mathrm{e}}^{-(t+ \mathrm{i}u)}({\mathrm{e}}^{-\mathrm{i}s}\varepsilon ) \big) = 
\big( 0, {\mathrm{e}}^{-t}({\mathrm{e}}^{-\mathrm{i}(s +u)}\varepsilon ) \big) . 
\end{displaymath}
So $\lim_{t \rightarrow \infty}{\varphi }^{\mathcal{H}}_{t + \mathrm{i}u}
(0, {\mathrm{e}}^{-\mathrm{i}s} \varepsilon ) = 
(0,0)$, that is, $(0, {\mathrm{e}}^{-\mathrm{i}s} \varepsilon ) \in W_s(0)$, since 
${\varphi }^{\mathcal{H}}_{t+ \mathrm{i}u }$ is ${\varphi }^{q_1}_t \comp {\varphi }^{q_2}_{-u}$ in real terms. Similarly,  $({\mathrm{e}}^{\mathrm{i}s} \varepsilon ,0) \in W_u(0)$. \hfill $\square $ \medskip 

\noindent \textbf{Lemma 3.5} With $c \ne 0$ the circle ${\mathcal{C}}_u = \{ {\eta }_s(c) \setrule \, 
s \in \R \} $ in ${\mathcal{H}}^{-1}(c)$ is a cross section for the flow 
${\varphi }^{\mathcal{H}}_{\tau }$ on ${\mathcal{H}}^{-1}(c)$. 
In other words, every complex integral curve of $X_{\mathcal{H}}$ starting at a point $(z^0_1, z^0_2) \in 
{\mathcal{H}}^{-1}(c)$ intersects ${\mathcal{C}}_u$. \medskip 

\noindent \textbf{Proof.} To show that the complex integral curve 
$\tau \mapsto {\varphi }^{\mathcal{H}}_{\tau}(z^0_1, z^0_2)$ intersects ${\mathcal{C}}_u$ we need to show 
that there is $( {\mathrm{e}}^{\mathrm{i}v}\varepsilon , \frac{c}{{\mathrm{e}}^{\mathrm{i}v}\varepsilon} ) \in {\mathcal{C}}_u$ and a $\tau = t+\mathrm{i}s$ such that 
\begin{displaymath}
( {\mathrm{e}}^{\mathrm{i}v}\varepsilon , \frac{c}{{\mathrm{e}}^{\mathrm{i}v}\varepsilon} ) = 
{\varphi }^{\mathcal{H}}_{t+ \mathrm{i}s}(z^0_1, z^0_2) = 
( {\mathrm{e}}^{t+ \mathrm{i}s}z^0_1 , {\mathrm{e}}^{-t -\mathrm{i}s}z^0_2 ). 
\end{displaymath}
It is enough to solve 
\begin{equation}
 {\mathrm{e}}^{t+ \mathrm{i}s}z^0_1 = {\mathrm{e}}^{\mathrm{i}v}\varepsilon , 
 \label{eq-s3three}
 \end{equation}
 because equation (\ref{eq-s3three}) implies ${\mathrm{e}}^{-t -\mathrm{i}s}z^0_2 = 
 \frac{c}{{\mathrm{e}}^{\mathrm{i}v}\varepsilon}$ since $z^0_1z^0_2 = c$ and $c \ne 0$. 
 Write $z^0_1 = |z^0_1| {\mathrm{e}}^{\mathrm{i}\, \mathrm{arg}\, z^0_1}$. Then (\ref{eq-s3three}) 
 reads 
 \begin{equation}
 {\mathrm{e}}^{t + \mathrm{i}(s + \mathrm{arg}\, z^0_1)}|z^0_1| = {\mathrm{e}}^{\mathrm{i}v}\varepsilon . 
 \label{eq-s3four}
 \end{equation}
 So $s = v - \mathrm{arg}\, z^0_1$ and $t = \ln \frac{\varepsilon }{|z^0_1|}$ solves (\ref{eq-s3four}). 
 Here $v$ may be chosen freely.  \hfill $\square $ \medskip  
 
\noindent \textbf{Lemma 3.6} Fix $c \ne 0$. Let ${\mathcal{C}}_s$ be the circle 
$\{ {\xi}_s(c) \setrule \, s \in \R \}$ in ${\mathcal{H}}^{-1}(c)$. The points ${\xi }_s(c) \in {\mathcal{C}}_s$ 
and ${\eta }_s(c) \in {\mathcal{C}}_u$ lie on a complex integral curve of $X_{\mathcal{H}}$. \medskip 

\noindent \textbf{Proof.} Set $\tau = - \ln \frac{c}{{\varepsilon }^2}$. Then 
\begin{align}
{\varphi }^{\mathcal{H}}_{\tau }\big( {\xi }_s(c) \big) & = 
\big( {\mathrm{e}}^{-\ln \frac{c}{{\varepsilon }^2}} (\ttfrac{c}{{\mathrm{e}}^{-\mathrm{i}s}\varepsilon } ) , 
{\mathrm{e}}^{\ln \frac{c}{{\varepsilon }^2}} ({\mathrm{e}}^{-\mathrm{i}s} \varepsilon )\big) 
= \big( {\mathrm{e}}^{\mathrm{i}s} \varepsilon , \ttfrac{c}{{\mathrm{e}}^{\mathrm{i}s}\varepsilon} \big)  
= {\eta }_s(c). \tag*{$\square $}
\end{align}

Since $\mathrm{Re}\, \tau =  2 \ln \varepsilon - \mathrm{Re} \ln c$, we get $\mathrm{Im}\, \tau = 
-\mathrm{Im} \ln c$. Thus $\mathrm{Im}\, \tau $ depends only on $c \ne 0$ and not on the 
choice of the circles ${\mathcal{C}}_{s,u}$ on ${\mathcal{H}}^{-1}(c)$ or on the choice of points 
${\xi }_s(c) \in {\mathcal{C}}_u$ or ${\eta }_s(c) \in {\mathcal{C}}_s$. Hence 
$\mathrm{Im}\, \tau $ is an intrinsic property of the flow of ${\varphi }^{\mathcal{H}}$ of the vector 
field $X_{\mathcal{H}}$ on ${\mathcal{H}}^{-1}(c)$. \medskip 

We look at the map 
\begin{displaymath}
\mu : {\C }^{\times } \rightarrow S^1: c \mapsto \mathrm{Im} \, \tau (c) = -\mathrm{Im}\ln c . 
\end{displaymath}

\noindent \textbf{Claim 3.7} The winding number of the map $\mu $ is $-1$, which is the \emph{scattering monodromy} of the standard focus-focus system. \medskip 

\noindent \textbf{Proof.} We now verify this last assertion. 
On ${\C}^{\times}$ with coordinate $z_1$ we have a real $1$-form 
$\vartheta = - \mathrm{Im}\frac{\dee z_1}{z_1} = - \mathrm{Im} \dee \, \ln z_1$. So  
\begin{align*}
\vartheta & =  - \mathrm{Im}\frac{\dee z_1}{z_1} 
= - \mathrm{Im} \Big( \frac{\dee x - \mathrm{i}\dee y}{x- \mathrm{i} y} \Big) \\
& = -(x^2+y^2)^{-1} \mathrm{Im} [ (x + \mathrm{i}y)(\dee x - \mathrm{i} \dee y)] \\
& = -(x^2+y^2)^{-1} (y \dee x - x \dee y ) =  \dee \, \big( {\tan }^{-1} \frac{x}{y} \big) . 
\end{align*}
Let $\pi : {\C}^{\times } \times \C \rightarrow {\C}^{\times }:
(z_1,z_2) \mapsto z_1$ be the projection map on the first factor. Set 
\begin{displaymath}
\theta = ({\pi }^{\ast }\vartheta )|{\mathcal{H}}^{-1}({\C}^{\times }). 
\end{displaymath}
Then $\theta $ is a closed $1$-form on ${\mathcal{H}}^{-1}({\C}^{\times })$ since 
$\dee \theta = \dee \, ({\pi }^{\ast }\vartheta )|{\mathcal{H}}^{-1}({\C}^{\times }) = 
\big( {\pi }^{\ast }(\dee \vartheta )\big) |{\mathcal{H}}^{-1}({\C}^{\times } )=0$. 
Because 
\begin{align*}
X_{\mathrm{Im}\mathcal{H}} \lefthook \, \theta & = X_{q_2} \lefthook \theta 
\\
& \hspace{-.65in}= \big( -y\frac{\partial }{\partial x} + x\frac{\partial }{\partial y} 
-p_y\frac{\partial }{\partial p_x} + p_x \frac{\partial}{\partial p_y} \big) 
\lefthook -(x^2+y^2)^{-1}(y \dee x - x \dee y ) \\
&\hspace{-.65in} = (x^2+y^2)^{-1}(y^2+x^2)= 1, 
\end{align*}
$\theta $ is a connection $1$-form on the bundle $\widehat{\rho }: {\mathcal{H}}^{-1}({\C}^{\times }) 
\mapsto {\C}^{\times }$. Since $-\mathrm{Im} \Big[ \frac{\dee {\mathrm{e}}^{\mathrm{i}s} z_1}
{{\mathrm{e}}^{\mathrm{i}s} z_1}\Big] $ $= - \mathrm{Im} \frac{\dee z_1}{z_1}$, the connection $1$-form 
$\theta $ is invariant under the flow of $X_{\mathrm{Im}\, \mathcal{H}}$ 
on ${\mathcal{H}}^{-1}({\C}^{\times })$. \medskip 

Let 
\begin{displaymath}
\Sigma = \widetilde{\Sigma }|S^1_r: S^1_r \subseteq 
{\C}^{\times } \rightarrow {\mathcal{H}}^{-1}(S^1_r): s \mapsto (s,1) 
\end{displaymath}
be a smooth section of the bundle $\widehat{\rho}$. Consider the real curve  
${\Gamma }_{\Sigma (s)}$ on ${\mathcal{H}}^{-1}(S^1_r)$, where 
${\Gamma }_{\Sigma (s)}(t) = {\varphi }^{\mathrm{Re}\, \mathcal{H}}_t\big( \Sigma (s) \big) $. 
With $(s,1) = {\Gamma }_{\Sigma (s)}(0)$ we have 
\begin{displaymath}
\big( z_1(t), z_2(t) \big) = 
{\Gamma }_{\Sigma (s)}(t) = {\varphi }^{q_1}_t(s,1) = ({\mathrm{e}}^t s, {\mathrm{e}}^t). 
\end{displaymath}
Thus for each $s \in S^1_r$ the scattering phase $\Theta (s)$ of ${\Gamma }_{\Sigma (s)}$ with 
respect to the connection $1$-form $\theta $ is $\int_{{\Gamma }_{\Sigma (s)}} \theta $. We have 
\begin{align*}
\Theta (s) & = \int_{{\Gamma }_{\Sigma (s)}} \langle \theta \! \mid 
\frac{\dee {\Gamma }_{\Sigma (s)}}{\dee t} \rangle \, \dee t = 
\int^{\infty}_{-\infty} \langle \theta \! \mid X_{q_1} \rangle \, \dee t \\
& = \int^{\infty}_{-\infty} \frac{\dee \theta }{\dee t}\big( {\Gamma }_{\Sigma (s)}(t) \big) \, \dee t , \, \, \, 
\parbox[t]{2.5in}{by definition of infinitesimal elevation, see \cite[p.436]{bates-cushman} } \\
& = - \int^{\infty}_{-\infty} \frac{\dee}{\dee t} \big( \mathrm{Im} \ln z_1(t) \big) \, \dee t = 
-\mathrm{Im} \ln z_1(t) \! \mid^{\infty}_{-\infty} \\ 
& = - \mathrm{Im} \ln ({\mathrm{e}}^t \, s) \! \mid^{\infty}_{-\infty} = 
-(\mathrm{Im} \ln s + \mathrm{Im} \ln {\mathrm{e}}^t) \! \mid^{\infty}_{-\infty} \\
& =  -\mathrm{Im} \ln s.
\end{align*}
This verifies the assertion. \hfill $\square $

\section{Connection $1$-form}

In this section we construct a connection $1$-form on ${\R }^4$, which is invariant under 
the flow of $X_K$. \medskip 

\noindent \textbf{Lemma 4.1} Let ${\gamma }_z$ be an integral curve of the 
vector field $X_{q_2} = -y \frac{\partial }{\partial x} + x \frac{\partial }{\partial y} 
-p_y \frac{\partial }{\partial p_x} + p_x \frac{\partial }{\partial p_y}$ on $\big( {\R }^4, \omega  = 
-\dee \alpha = -\dee \, ( p_x\dee x + p_y \dee y) \big) $ starting at $z = (x,y,p_x,p_y)$. Then 
\begin{equation}
q_2(z) = \frac{1}{2\pi } \int_{{\gamma }_z} \alpha .
\label{eq-s4one}
\end{equation}

\noindent \textbf{Proof.} We compute. By definition 
\begin{align*}
\frac{1}{2\pi } \int_{{\gamma }_z} \alpha & = 
\frac{1}{2\pi } \int^{2\pi}_0 \langle \alpha \! \mid \frac{\dee {\varphi }^{q_2}_s}{\dee s} \rangle (z) \dee s 
= \frac{1}{2\pi } \int^{2\pi }_0 \langle \alpha \! \mid X_{q_2} \rangle \big( {\varphi }^{q_2}_s(z) \big) \dee s .
\end{align*}
But $\langle \alpha \! \mid X_{q_2} \rangle = X_{q_2} \lefthook \alpha = 
-yp_x +xp_y = q_2$. So 
\begin{displaymath}
\frac{1}{2\pi } \int_{{\gamma }_z} \alpha = \frac{1}{2\pi } \int^{2\pi }_0 q_2\big( {\varphi }^{q_2}_s(z) \big) \dee s 
= q_2(z), 
\end{displaymath}
since $q_2 $ is an integral of $X_{q_2}$. \hfill $\square $ \medskip 

Let ${\Phi }_t = {\iota }^{-1} \comp {\varphi }^X_t \comp \iota $, where ${\varphi }^X_t$ for $t \in [0,1]$ is 
the flow of the vector field $X$ constructed in the proof of the complex Morse lemma for the 
complex function $\mathcal{H}$ on $[0,1] \times {\C}^2$. \medskip 

\noindent \textbf{Lemma 4.2} For $t \in [0,1]$ let  
\begin{displaymath}
K_t(z) = \frac{1}{2\pi } \int_{{\gamma }_z} {\alpha}_t = 
\frac{1}{2\pi } \int_{{\gamma }_z}({\Phi }^{-1}_t)^{\ast }\alpha . 
\end{displaymath}
Then $K_t = ({\Phi }^{-1}_t)^{\ast }q_2$. Since ${\Phi }_1 = \Phi $, it follows that $K_1 =K$. \medskip 

\noindent \textbf{Proof.} We compute. From lemma 4.1 we obtain 
\begin{align}
({\Phi }^{-1}_t)^{\ast }q_2(z) & = \frac{1}{2\pi }\int_{{\Phi }^{-1}_t \comp {\gamma }_z} \alpha 
= \frac{1}{2\pi} \int_{{\gamma }_z} ({\Phi }^{-1}_t)^{\ast }\alpha 
= \frac{1}{2\pi }\int_{{\gamma }_z} {\alpha }_t = K_t(z). \tag*{$\square $}
\end{align}

\noindent \textbf{Lemma 4.3} For $t \in [0,1]$ let ${\omega }_t = ({\Phi }^{-1}_t)^{\ast } \omega $. 
Since ${\Phi }_1 = \Phi $, it follows that ${\omega }_1 = ({\Phi }^{-1})^{\ast }\omega = \widetilde{\omega }$. Then 
\begin{equation}
X_{K_t} = ({\Phi }^{-1}_t)^{\ast }X_{q_2}. 
\label{eq-s4two}
\end{equation}

\noindent \textbf{Proof.} We compute. By definition 
\begin{align*}
X_{K_t} \lefthook {\omega }_t & = \dee K_t = \dee \, \big( ({\Phi }^{-1}_t)^{\ast }q_2 \big), \, \, \, 
\mbox{by lemma 4.1} \\
& = ({\Phi }^{-1}_t)^{\ast }\dee q_2 = ({\Phi }^{-1}_t)^{\ast }( X_{q_2} \lefthook \omega ) \\
& = ({\Phi }^{-1}_t)^{\ast }X_{q_2} \lefthook ({\Phi }^{-1}_t)^{\ast } \omega = 
({\Phi }^{-1}_t)^{\ast }X_{q_2} \lefthook {\omega }_t. 
\end{align*}
So $X_{K_t} = ({\Phi }^{-1}_t)^{\ast }X_{q_2}$, since ${\omega }_t$ is nondegenerate. 
\hfill $\square $ \medskip 

Let $\theta = {\pi }^{\ast }\big( \dee \, {\tan }^{-1} \frac{x}{y} \big) $. Then $\theta $ is a connection 
$1$-form on $({\R }^4, \omega )$, because $\theta $ is closed, $X_{q_2} \lefthook \theta = 1$ and 
$\theta $ is ${\varphi }^{q_2}_s$-invariant. \medskip 

\noindent \textbf{Claim 4.4} For $t \in [0,1]$ let ${\theta }_t = ({\Phi }^{-1}_t)^{\ast } \theta $. Then 
${\theta }_t$ is a connection $1$-form on $({\R }^4, {\omega }_t)$. In particular, since ${\Phi }_0 = 
\mathrm{id}$, we obtain ${\theta }_0 = \theta $, which is a connection $1$-form on $({\R }, \omega )$. Also
${\theta }_1 = ({\Phi }^{-1})^{\ast} \theta $ is a connection $1$-form on $({\R }^4, \widetilde{\omega })$. 
\medskip 

\noindent \textbf{Proof.} By definition 
\begin{align*}
{\theta }_t & = ({\Phi }^{-1}_t)^{\ast } \theta  = 
({\Phi }^{-1}_t)^{\ast } {\pi }^{\ast } \big( \dee \, {\tan}^{-1}\frac{x}{y} \big) = 
\dee \big( ({\Phi }^{-1}_t)^{\ast } {\pi }^{\ast } \big( {\tan}^{-1}\frac{x}{y} \big) . 
\end{align*}
So ${\theta }_t$ is a closed $1$-form. Next $X_{K_t} \lefthook {\theta }_t =1$, because 
\begin{align*}
X_{K_t} \lefthook {\theta }_t & = \langle {\theta }_t \! \mid X_{K_t} \rangle = 
\langle ({\Phi }^{-1}_t)^{\ast} \theta \! \mid ({\Phi }^{-1}_t)^{\ast }X_{q_2} \rangle \\
& = ({\Phi }^{-1}_t)^{\ast } \big( \langle \theta \! \mid X_{q_2} \rangle \big) = 
({\Phi }^{-1}_t)^{\ast } ( X_{q_2} \lefthook \theta ) \\
& = ({\Phi }^{-1}_t)^{\ast } 1 = 1.
\end{align*}
Let ${\varphi }^{K_t}_s$ be the flow of $X_{K_t}$. Since $X_{K_t} = ({\Phi }^{-1}_t)^{\ast }X_{q_2}$, we obtain  
${\varphi }^{K_t}_s = {\Phi }_t \comp {\varphi }^{q_2}_s \comp {\Phi }^{-1}_t$. The $1$-form 
${\theta }_t$ is invariant under ${\varphi }^{K_t}_s$. To see this we compute 
\begin{align*}
({\varphi }^{K_t}_s)^{\ast }{\theta }_t &= 
({\Phi }_t \comp {\varphi }^{q_2}_s \comp {\Phi }^{-1}_t)^{\ast }({\Phi }^{-1}_t)^{\ast } \theta  = 
({\Phi }^{-1}_t)^{\ast } \big( ({\varphi }^{q_2}_s)^{\ast} \theta \big) = 
({\Phi }^{-1}_t)^{\ast } \theta = {\theta }_t.
\end{align*}
Thus ${\theta }_t$ is a connection $1$-form on $({\R }^4, {\omega }_t)$. \hfill $\square $ \medskip

We are now in position to prove the geometric scattering monodromy theorem. 
For each $t \in [0,1]$ consider the integrable Hamiltonian system 
$(q_1, K_t, {\R }^4, {\omega }_t)$, where $K_t = ({\Phi }^{-1}_t)^{\ast }q_2$. Let 
$S^1_r$ be a circle of radius $r$ in $U \setminus \{ (0,0) \} \subseteq {\R }^2$. Consider the bundle 
\begin{displaymath}
{\widehat{\rho}}_t = F_t| F^{-1}_t(S^1_r): F^{-1}_t(S^1_r) \subseteq W \rightarrow S^1_r \subseteq 
U \setminus \{ (0,0) \} : z \mapsto \big( q_1(z), K_t(z) \big) , 
\end{displaymath}
which has a global section 
\begin{displaymath}
{\Sigma }_t: S^1_r \subseteq {\R }^2\rightarrow F^{-1}_t(S^1_r): s \mapsto 
({\Phi }_t \comp \widehat{\Sigma })(s),  
\end{displaymath}
where 
\begin{equation}
\widehat{\Sigma } = {\iota }^{-1}\comp \Sigma \comp j: S^1_r \subseteq {\R }^2 
\rightarrow F^{-1}_t(S^1_r).
\label{eq-s4twostar} 
\end{equation}
For each $s\in S^1_r$ consider the curve ${\Gamma }_{{\Sigma }_t(s)}$ on $F^{-1}_t(S^1_r)$ given by 
${\Gamma }_{{\Sigma }_t(s)}(v) = {\varphi }^{K_t}_v\big( {\Sigma }_t(s) \big)$. \medskip

\noindent \textbf{Claim 4.5} For each $t \in [0,1]$ the scattering phase 
${\Theta }_t(s)$ of the curve ${\Gamma }_{{\Sigma }_t(s)}$ on $F^{-1}_t(S^1_r)$ 
with respect to the connection $1$-form ${\theta }_t$ is $\Theta (s)$. \medskip 

\noindent \textbf{Proof.} We compute. 
\begin{align*}
{\Theta }_t(s) & = \int_{{\Gamma }_{{\Sigma }_t(s)}} {\theta }_t = 
\int_{{\Gamma }_{{\Sigma }_t(s)}} ({\Phi }^{-1}_t)^{\ast } \theta 
= \int_{{\Phi }^{-1}_t \comp {\Gamma }_{{\Sigma }_t(s)} } \theta = 
\int_{{\Gamma }_{\widehat{\Sigma }(s)}} \theta 
= \Theta (s). 
\end{align*}
The second to last equality above follows because ${\Gamma }_{\widehat{\Sigma }(s)}(v) = 
{\varphi }^{q_2}_v\big( \widehat{\Sigma }(s) \big) $, where $\widehat{\Sigma }$ (\ref{eq-s4twostar}) 
is a global section of the bundle ${\widehat{\rho}}_t$, gives 
\begin{align}
{\Gamma}_{{\Sigma }_t(s)}(v) & = {\varphi }^{K_t}_v \big( {\Sigma }_t(s) \big) \notag \\
& = ({\Phi }_t \comp {\varphi }^{q_2}_v \comp {\Phi }^{-1}_t)\big( {\Sigma }_t(s) \big) , \, \, \,  
\mbox{because $X_{K_t} = ({\Phi }^{-1}_t)^{\ast }X_{q_2}$}  \notag \\
& = {\Phi }_t \big( {\varphi }^{q_2}_v \big( \widehat{\Sigma }(s) \big), \, \, \, 
\mbox{since ${\Sigma }_t = {\Phi }_t \comp \widehat{\Sigma }$} \notag  \\
& = {\Phi }_t\big( {\Gamma }_{\widehat{\Sigma }(s)}(v) \big) .  \tag*{$\square $}
\end{align}

\noindent \textbf{Corollary 4.5A} For each $t \in [0,1]$ the mapping 
${\Theta }_t: S^1_r \rightarrow \C : s \mapsto {\Theta }_t(s) $ has 
winding number $-1$. \medskip 

\noindent \textbf{Proof.} This follows because the scattering map $\Theta $ of 
the standard focus-focus system has winding number $-1$. \hfill $\square $ \medskip 

\vspace{-.15in}We now complete the proof of the geometric scattering monodromy theorem. The winding 
number of the scattering map ${\Theta }_1$ with respect to the connection $1$-form 
${\theta }_1$ on $({\R }^4, {\omega }_1 = \widetilde{\omega })$ associated to the 
family of curves ${\Gamma }_{{\Sigma }_1(s)}$ on $F^{-1}_1(S^1_r)$ is $-1$. This proves 
the geometric scattering monodromy theorem, becuase $X_I = {\Phi }^{\ast }X_K$. So 
the scattering phase map associated to the scattering phase of 
${\Phi }^{-1}\comp {\Gamma }_{\widehat{\Sigma }(s)}$ with respect to the connection $1$-form 
${\Phi }^{\ast }{\theta }$ has the same winding number 
as the scattering phase map associated to the scattering phase of 
${\Gamma }_{\widehat{\Sigma }(s)}$ with respect to the connection $1$-form $\theta $, since 
$\Phi $ is a diffeomorphism.

\section{Appendix}

In this appendix we prove a complex version of the Morse lemma. Our proof was inspired by 
the proof of the focus-focus Morse lemma in \cite{vungoc-wacheux}. \medskip

\noindent \textbf{Lemma A1} (Morse lemma). Let 
\begin{displaymath}
\mathcal{H}: [0,2] \times {\C }^2 \rightarrow \C: (t,z) \mapsto Q(z) + t R(z), 
\end{displaymath}
where $Q$ is a nondegenerate homogeneous quadratic polynomial and $R$ is a smooth 
function, which is flat to second order at the origin $(0,0)$. Then there is an open 
neighborhood $U$ of $(0,0)$ in ${\C}^2$ and a diffeomorphism $\Phi $ of $U$ into itself 
with $\Phi (0,0) = (0,0)$ such that ${\Phi }^{\ast }{\mathcal{H}}_1 = Q$ on $U$. Moreover, 
$\Phi $ is isotopic to ${\mathrm{id}}_U$. \medskip 

\noindent \textbf{Proof.} By a complex linear change of coordinates we may assume that 
$Q(z) = z_1z_2$. We want to find a time dependent vector field $X = X_t + \frac{\partial }{\partial t}$ on 
$(0,2) \times {\C}^2$ whose flow ${\varphi }^X_t$ satisfies 
\begin{equation}
({\varphi }^X_t)^{\ast }\mathcal{H} = Q, \, \, \, 
\mbox{for every $t \in [0,1]$.}
\label{eq-Aone}
\end{equation}
Differentiating (\ref{eq-Aone}) gives $0 = ({\varphi }^X_t)^{\ast }\big( \frac{\partial \mathcal{H}}{\partial t} 
+ L_{X_t}{\mathcal{H}}_t \big)$. Since $\frac{\partial \mathcal{H}}{\partial t} = R$, we need to find a vector field $X_t$ on ${\C}^2$ such that 
\begin{equation}
\dee {\mathcal{H}}_t(z) X_t(z) = - R(z) \, \, \, \mbox{for all $t \in [0,1]$.}
\label{eq-Atwo}
\end{equation}

Now 
\begin{displaymath}
\dee H_t(z) = \big( z_2 + t \frac{\partial R}{\partial z_1}(z) \big) \dee z_1 + 
\big( z_1 + t \frac{\partial R}{\partial z_2}(z) \big) \dee z_2. 
\end{displaymath}
For some smooth functions $A$ and $B$ on $(0,2) \times {\C}^2$
\begin{equation}
X_t(z) = A(t,z) \frac{\partial }{\partial z_1} + B(t,z) \frac{\partial }{\partial z_2}.
\label{eq-Athree}
\end{equation}
Since $R(0) =0$, by the integral form of Taylor's theorem $R(z) = G_1(z) z_1 + G_2(z) z_2$, 
where $G_j$ are smooth functions with $G_j(0) = \frac{\partial R}{\partial z_j}(0)$ for $j=1,2$. Since 
$R$ is flat to second order at $0$, we get $G_j(0) =0$. Thus (\ref{eq-Atwo}) can be written as 
\begin{align}
-G_1(z) z_1 - G_2(z) z_2 & = A(t,z) \big( z_2 + t\frac{\partial R}{\partial z_1}(z) \big) + 
B(t,z) \big( z_1 + \frac{\partial R}{\partial z_2}(z) \big) . 
\label{eq-Afour}
\end{align} 
Again by Taylor's theorem, $\frac{\partial R}{\partial z_j}(z) = F_j(z) z_1 + E_j(z) z_2$ for $j =1,2$, 
where $F_j(0) = \frac{{\partial }^2R}{\partial z_1 \, \partial z_j}(0)$ and $E_j(0) = 
\frac{{\partial }^2R}{\partial z_2 \, \partial z_j}(0)$. Since $R$ is flat to second order at $0$, 
it follows that $F_j(0) =0$ and $E_j(0) =0$. Thus equation (\ref{eq-Afour}) becomes 
\begin{align*}
-G_1(z) z_1 - G_2(z) z_2 & = A(t,z) \big( z_2 +t(F_1(z) z_1 + E_1(z) z_2) \big) \\
& \hspace{.25in}+ B(t,z) \big( z_1 + t (F_2(z)z_1 + E_2(z) z_2) \big)  .
\end{align*}
Equating the coefficients of $z_1$ and $z_2$ in the equation above, we get  
\begin{displaymath}
-\mbox{\footnotesize $\begin{pmatrix} G_1(z) \\ G_2(z) \end{pmatrix}$} = 
\mbox{\footnotesize $\begin{pmatrix} t F_1(z) & 1+t F_2(z) \\ 1+t E_1(z) & tE_2(z) \end{pmatrix} \, 
\begin{pmatrix} A(t,z) \\ B(t,z) \end{pmatrix} $} = \mathcal{A}(t, z)\mbox{\footnotesize $\begin{pmatrix} 
A(t,z) \\ B(t,z) \end{pmatrix}$.}
\end{displaymath}
So 
\begin{align*}
|\det \mathcal{A}(t,z)| & = |1 +t\big( E_1(z) + F_2(z) \big) + t^2 \big( E_1(z)F_2(z)-E_2(z)F_2(z) \big) | \\
& \ge 1 - |t| \, |\big( E_1(z) + F_2(z) \big) + t^2 \big( E_1(z)F_2(z)-E_2(z)F_2(z) \big) | . 
\end{align*}
Let $U$ be an open neighborhood of $0 \in {\C}^2$ such that for $i=1,2$ 
\begin{equation}
|E_i(z)| < \ttfrac{1}{16} \, \, \, \mathrm{and} \, \, \, |F_i(z)| < \ttfrac{1}{16}.
\label{eq-Afive}
\end{equation}
Then 
\begin{align*}
|t| \, |\big( E_1(z) + F_2(z) \big) + t^2 \big( E_1(z)F_2(z)-E_2(z)F_2(z) \big) | &  \\
&\hspace{-2.5in} \le |t| \big( |E_1(z)| + |F_2(z)| + |t| \big[ |E_1(z)| \, |F_2(z)| + |E_2(z)|\, |F_1(z)| \big] \big) \\ 
& \hspace{-2.5in} < 2\big[ \ttfrac{1}{16} + \ttfrac{1}{16} + 2(\ttfrac{1}{16}\cdot \ttfrac{1}{16} + 
\ttfrac{1}{16}\cdot \ttfrac{1}{16}) \big], \, \, \mbox{using (\ref{eq-Afive}) and $t\in [0,2]$} \\
& \hspace{-2.5in} = \ttfrac{17}{64} . 
\end{align*}
Thus the matrix $\mathcal{A}(t,z)$ is invertible for all $z \in U$ and all $t \in [0,2]$. With {\tiny \raisebox{2pt}{$\begin{pmatrix} A(t,z) \\ B(t,z) \end{pmatrix}$}}$= -{\mathcal{A}(t,z)}^{-1}${\tiny \raisebox{1pt}{$\begin{pmatrix} G_1(z) \\ G_2(z) \end{pmatrix}$}} we have determined the 
vector field $X_t$ (\ref{eq-Athree}) on $(0,2) \times {\C}^2$ which solves equation (\ref{eq-Atwo}). \medskip 

Because $X_t(0,0) = (0,0)$ we can shrink $U$ if necessary so that the flow ${\varphi }^X_t$ of the 
vector field $X$ sends $[0,1] \times U$ to $U$. Set $\Phi = {\varphi }^X_1$. Then 
${\Phi }^{\ast }{\mathcal{H}}_1 =({\varphi }^X_1)^{\ast }{\mathcal{H}}_t = Q$ on $U$. The 
diffeomorphism $\Phi $ is isotopic to ${\mathrm{id}}_U$, since it is the time $1$ map 
of a flow. \hfill $\square $ \bigskip

\noindent \textbf{Statements and Declarations}
\par \noindent \hspace{.25in}\parbox[t]{4.5in}{The author received no funding for this research. This research involved no animals or humans and neither generated or used any computer programs or data.}

\end{document}